\documentclass[12pt]{amsart}
 
\input{xy} 
\xyoption{all}
\usepackage{graphicx}
\usepackage{color}
\usepackage{amssymb}

\newtheorem{theorem}{Theorem} [section]
\newtheorem{prop}[theorem]{Proposition} 
\newtheorem{lemma}[theorem]{Lemma}
\newtheorem{cor}[theorem]{Corollary}

\numberwithin{equation}{section} 
\numberwithin{figure}{section}

\newcommand\C{{\mathbb C}}

\renewcommand\P{{\mathbb P}}
\newcommand\R{{\mathbb R}}
\newcommand\Z{{\mathbb Z}}

\newcommand\N{{\mathbb N}}

\newcommand\Hyp{{\mathbb H}}
\newcommand\cM{\mathcal{M}}
\newcommand\cB{\mathcal{B}}
\newcommand\cC{\mathcal{C}}
\newcommand{\cT}{{\mathcal T}}

\newcommand{\cG}{\mathcal {G}}
\newcommand{\cH}{\mathcal {H}}

\newcommand{\cP}{\mathcal{P}}
\newcommand{\cS}{\mathcal{S}}
\newcommand{\cW}{\mathcal{W}}
\newcommand\eps{\varepsilon}

\renewcommand\phi{\varphi}

\newcommand\GL{\mathrm{GL}}

\newcommand\iso{\simeq} 

\renewcommand\mod{\operatorname{mod}}  

\newcommand\rank {\operatorname{rank}}
\newcommand\del{\partial} 
\newcommand\delbar{\bar{\del}} 

\renewcommand\Im {\operatorname{Im}} 
 
\newcommand\MP {\mathcal{M}}

\newcommand{\mtwo}[4]                            
{\mbox{$\left(\begin{array}{cc}                  
#1 & #2 \\
#3 & #4 
\end{array}
\right)$}}

\begin{document}

\title{Critical heights on the moduli space of polynomials}

\author{Laura DeMarco and Kevin Pilgrim}

\begin{abstract}  
Let $\MP_d$ be the moduli space of one-dimensional complex polynomial dynamical systems.  The escape rates of the critical points determine a {\em critical heights map} $\cG: \MP_d \to \R^{d-1}$.  For generic values of $\cG$, each connected component of a fiber of $\cG$ is the deformation space for twist deformations on the basin of infinity.  We analyze the quotient space $\cT_d^*$ obtained by collapsing each connected component of a fiber of $\cG$ to a point.  The space $\cT_d^*$ is a parameter-space analog of the polynomial tree $T(f)$ associated to a polynomial $f:\C\to\C$, studied in \cite{DM:trees}, and there is a natural projection from $\cT_d^*$ to the space of trees $\cT_d$.  We show that the projectivization $\P\cT_d^*$ is compact and contractible; further, the shift locus in $\P\cT_d^*$ has a canonical locally finite simplicial structure.  The top-dimensional simplices are in one-to-one corespondence with topological conjugacy classes of structurally stable polynomials in the shift locus.  
\end{abstract}

\date{\today}

\maketitle


\thispagestyle{empty}



\section{Introduction}

This article continues the study, initiated by Branner and Hubbard in \cite{Branner:Hubbard:1, Branner:Hubbard:2},  of the global structure of the moduli space of polynomials $f:\C\to\C$ determined by the dynamics on the basin of infinity,
	$$ X(f) = \{z\in\C: f^n(z) \to \infty \mbox{ as } n\to \infty\}.$$
Associated to a polynomial $f$ are certain analytic invariants, such as the coordinates of its critical points in the uniformizing B\"ottcher coordinates on $X(f)$, and certain combinatorial invariants, such as the simplicial tree $T(f)$ defined in \cite{DM:trees} or the tableau of \cite{Branner:Hubbard:2}.   In this article and its sequel \cite{DP:combinatorics}, we investigate the extent to which these invariants determine the dynamics of the polynomial $f$ on its basin of infinity.  

This study is partly motivated by a desire to understand the {\em shift locus} in the moduli space, the set of polynomials $f$ where all critical points lie in $X(f)$.  In the shift locus, the conformal conjugacy class of a polynomial is uniquely determined by its restriction to the basin of infinity, and all such polynomials are topologically conjugate on their Julia sets.  The shift locus has complicated topology, which can be seen in the rich variation of the global dynamics of these polynomials \cite{Blanchard:Devaney:Keen}, and there is currently no known invariant which classifies its topological conjugacy classes.  In a sequel, we provide a combinatorial method to complete the classification of polynomial dynamical systems restricted to their basins of infinity.   Here we concentrate on the general topological features of conjugacy classes and how they fit together in the moduli space.

\subsection{Critical heights}
Let $\MP_d$ denote the space of conformal conjugacy classes of polynomials $f:\C\to\C$ of degree $d>1$. 
We consider the {\em critical heights map} 
	$$\mathcal{G}: \MP_d \to \R^{d-1}$$
defined by $\mathcal{G}(f) = (G_f(c_1), \ldots, G_f(c_{d-1}))$, where 
	$$G_f(z) = \lim_{n\to\infty} \frac{1}{d^n} \log^+ |f^n(z)|$$
is the escape-rate function of $f$, and $\{c_1, \ldots, c_{d-1}\}$ are the critical points of $f$, ordered so that $G_f(c_i)\geq G_f(c_{i+1})$.   Note that the connectedness locus $\mathcal{C}_d$, consisting of all polynomials with connected Julia set, is the fiber of $\cG$ over the origin in $\R^{d-1}$ .

Consider the subset 
	$$\cH_d = \{ (h_1, \ldots, h_{d-1}): h_1 \geq h_2 \geq \ldots \geq h_{d-1} \geq 0\}$$ 
of $\R^{d-1}$.  We begin with the following observation:

\begin{theorem}  \label{thm:surjective}
The critical heights map 
	$$\cG: \MP_d \to \cH_d$$ 
is continuous, proper and surjective.
\end{theorem}

\noindent
The continuity and properness of $\cG$ are well-known, following easily from \cite{Branner:Hubbard:1}.  Surjectivity follows from basic properties of analytic maps,  as $\cG$ is pluriharmonic where the critical points $c_i$ can be consistently labeled and have positive heights.  Details are given in \S\ref{polynomials}.

We study the decomposition of $\MP_d$ into connected components of the fibers of $\mathcal{G}$.   This critical-heights decomposition is related to twisting deformations on the basin of infinity:

\begin{theorem}  \label{generic fibers}
For generic values of $\cG: \MP_d \to \cH_d$, the fiber of $\cG$ is a finite, disjoint union of smooth $(d-1)$-dimensional tori; each connected component coincides with a twist-deformation orbit.  In particular, the connected components of a generic fiber are precisely the topological conjugacy classes of polynomials within that fiber.  
\end{theorem}

The generic values of Theorem \ref{generic fibers} have a precise characterization:  we show that a value $(h_1, \ldots, h_{d-1})$ of $\cG$ satisfies the conclusion of Theorem \ref{generic fibers} if and only if $h_i>0$ for all $i$ and $h_i \not= d^n h_j$ for all $i\not=j$ and $n\in \Z$.  Dynamically, these generic critical heights correspond to polynomials in the shift locus with all critical points in distinct foliated equivalence classes.  The {\em twist deformation}, introduced in \cite{McS:QCIII}, is a quasiconformal deformation of the basin of infinity of a polynomial which twists its fundamental annuli and preserves its critical heights; the turning curves of \cite{Branner:Hubbard:1} are a special case.  See \S\ref{sec:qc}.

\subsection{Monotone-light factorization}  A closed, proper map between Hausdorff spaces factors canonically as a {\em monotone} map (i.e. one whose fibers are connected) followed by a {\em light} map (one whose fibers are totally disconnected); see  \cite[Theorem I.4.3]{daverman:decompositions}.   Below, we apply this to the critical heights map.  

The critical heights map factors as a composition of continuous, proper, and surjective maps
 	$$\MP_d \to  \mathcal{B}_d \to \cT_d \to \cH_d$$
where $\mathcal{B}_d$ is the space of conformal conjugacy classes of polynomials restricted to their basins of infinity, introduced in \cite{DP:basins}, and $\cT_d$ is the space of polynomial trees introduced in \cite{DM:trees}.  We remark that this space $\cT_d$ differs somewhat from its original definition, in that we have included the unique {\em trivial tree} associated to polynomials with connected Julia set.  See \S\ref{sec:trees} for details.

The main result of \cite{DP:basins} states that the projection $\MP_d \to  \mathcal{B}_d$ is monotone, meaning that its fibers are connected.  On the other hand, the critical heights map on the space of trees $\cT_d \to \cH_d$ is light (Lemma \ref{THfibers}).  

The fibers of the projection $\mathcal{B}_d \to \cT_d$ can be disconnected, and the reason is not too surprising.  In forming the tree $T(f)$ of a polynomial $f$, connected components of the level sets of $G_f$ are collapsed to points, thus ``forgetting'' the complicated planar graph structure of the singular components.  We give explicit examples in the sequel \cite{DP:combinatorics}, but we remark here that even for generic height values in degree $d=3$, there are distinct topological conjugacy classes of polynomials with the same tree.

We define $\cT^*_d$ to be the quotient space of $\MP_d$ obtained by collapsing connected components of the fibers of $\cG$ to points; it is a parameter-space analog of the tree for a single polynomial.  The fact that $\cT_d \to \cH_d$ is light  implies that we obtain a new sequence of continuous, proper, surjective maps 
\begin{equation} \label{space sequence}
  \MP_d \to  \mathcal{B}_d \to \cT_d^* \to \cT_d \to \cH_d
\end{equation}
through which the critical heights map $\cG$ must factor.  By construction, the quotient space $\cT_d^*$ is the unique monotone-light factor for $\cG$.  That is, the fibers of $\MP_d\to\cT_d^*$ are connected (and, in fact, always contain a non-degenerate continuum), while the fibers of $\cT_d^*\to \cH_d$ are totally disconnected.  Moreover:

\begin{theorem}   \label{thm:TstarTfibers}
The fibers of the projection $\cT^*_d \to \cT_d$ are totally disconnected; the fibers are finite over the shift locus in $\cT_d$.  Furthermore, for each $h\geq 0$, the fiber over the unique tree with uniform critical heights $(h, h, \dots, h)$ is a single point.
\end{theorem}

From the definitions we deduce:

\begin{cor}
The space $\cT_d^*$ is the unique monotone-light factor for the tree projection map $\MP_d \to \cT_d$.  
\end{cor}

The objects of the space $\cT^*_d$ are at present somewhat mysterious; the notation $\cT_d^*$ has been chosen because elements should represent polynomial trees, augmented by a certain amount of combinatorial information which is finite (though unbounded) over the shift locus.  The extra information needed to specify an element of $\cT^*_d$ from the corresponding element of $\cT_d$ is the subject of \cite{DP:combinatorics}.  In particular, we provide examples which show that the projection $\cT^*_d \to \cT_d$ is a bijection if and only if $d=2$.   Theorem \ref{generic fibers} shows that $\cT_d^*$ is generically the orbit space for the twisting deformation; in \cite{DP:hausdorff} we showed that the full space $\cT_d^*$ can be interpreted as the {\em hausdorffization} of the twist-orbit space.

\subsection{Cone structures and stretching}
The Branner-Hubbard stretching operation, while not continuous on all of $\MP_d$, nevertheless descends to a well-defined, continuous action of $\R_+$ on each of the quotient spaces in (\ref{space sequence}), which is free and proper on the complement of the connectedness locus (Lemma \ref{lemma:stretching_continuous}).   The $\R_+$-actions are equivariant with respect to the projection maps in (\ref{space sequence}).  The heights map $\cG$ satisfies $$\cG(s\cdot f) = s\, \cG(f)$$ for all $s \in \R_+$.  The connectedness locus of each quotient space is the unique fixed point of the stretching action, and it is the limit of every orbit as $s$ decreases to $0$.  We obtain:

\begin{theorem}
\label{thm:cone}
The stretching operation induces a cone structure on each of the spaces $\mathcal{B}_d$, $\cT_d^*$, $\cT_d$, and $\cH_d$.  Specifically, they are cones over the compact sets  $\mathcal{B}_d(1)$, $\cT_d^*(1)$, $\cT_d(1)$, and $\cH_d(1)$, consisting of points with maximal critical escape rate equal to $1$, with origin at the connectedness locus.    
\end{theorem}

\noindent
Because of the cone structure, there is a natural {\em projectivization} of each quotient space of $\MP_d$ in (\ref{space sequence}), identified with the slice with maximal critical height 1.  For the space of trees $\cT_d$, this projectivization is the space $\P\cT_d$ of \cite{DM:trees}.

In \cite{DP:basins}, we introduced a deformation of polynomials which ``flows" critical points upward along external rays until they reach the height of the fastest escaping critical point.  While this cannot be canonically defined on $\MP_d(1)$ or $\mathcal{B}_d(1)$ (even in the shift locus, due to the collision of external rays), the ambiguity is avoided when passing to the quotient space $\cT_d^*$. Using this upward flow, we deduce:

\begin{theorem}  \label{thm:contractible}
The projectivization $\P\cT_d^*$ is compact and contractible.
\end{theorem}

\noindent
The corresponding statement for the projectivized space of trees $\P\cT_d$ was proved in \cite{DM:trees}.

\subsection{Stratification of the shift locus}  \label{strat}
The shift locus in each of the spaces of (\ref{space sequence}) is the subset of points with all critical heights positive.  It forms a dense open subset of each of these quotient spaces of $\MP_d$.  Restricting to the shift locus yields a sequence of surjective proper maps 
\begin{equation}
\label{eqn:shift_chain}
 \cS\MP_d \to \cS\mathcal{B}_d \to \cS\cT_d^* \to \cS\cT_d \to \cS\cH_d.
 \end{equation}
The first map is a homeomorphism \cite[Theorem 1.1]{DP:basins},  and $\cS\cH_d$ is the subset of $\cH_d$ with all coordinates positive.

Each of the spaces in (\ref{eqn:shift_chain}) is stratified by the maximal number of independent critical heights, where positive real numbers $x$ and $y$ are {\em independent in degree $d$} if there is no integer $n$ such that $x = d^n y$.  Dynamically, the $N$-th stratum $\cS\MP_d^{N}$ consists of polynomials with exactly $N$ distinct critical foliated equivalence classes.  By the quasiconformal deformation theory of \cite{McS:QCIII}, the connected components of the open stratum $\cS\MP_d^{d-1}$ are precisely the topological conjugacy classes of structurally stable polynomials in the shift locus.

We denote the $N$-th stratum on each of the spaces in (\ref{eqn:shift_chain}) with the superscript $N$.  
The strata $\cS\cT_d^{*N}$, $\cS\cT_d^{N}$, and $\cS\cH_d^{N}$ are non-connected real manifolds, with each component homeomorphic to $\R^{N-1}$. The maps in (\ref{eqn:shift_chain}) preserve the strata, and the strata are invariant under stretching.  Passing to the quotients by stretching, which we denote with the prefix $\P$ for projectivization, we obtain a new sequence of surjective proper maps 
\begin{equation}
\label{eqn:Pshift_chain}
 \P\cS\MP_d \to \P\cS\mathcal{B}_d \to \P\cS\cT_d^* \to \P\cS\cT_d \to \P\cS\cH_d. 
 \end{equation}
It follows that the stratifications descend to the projectivized shift loci, and that the maps in (\ref{eqn:Pshift_chain}) also preserve the corresponding strata.

\begin{figure} 
\includegraphics[width=5.5in]{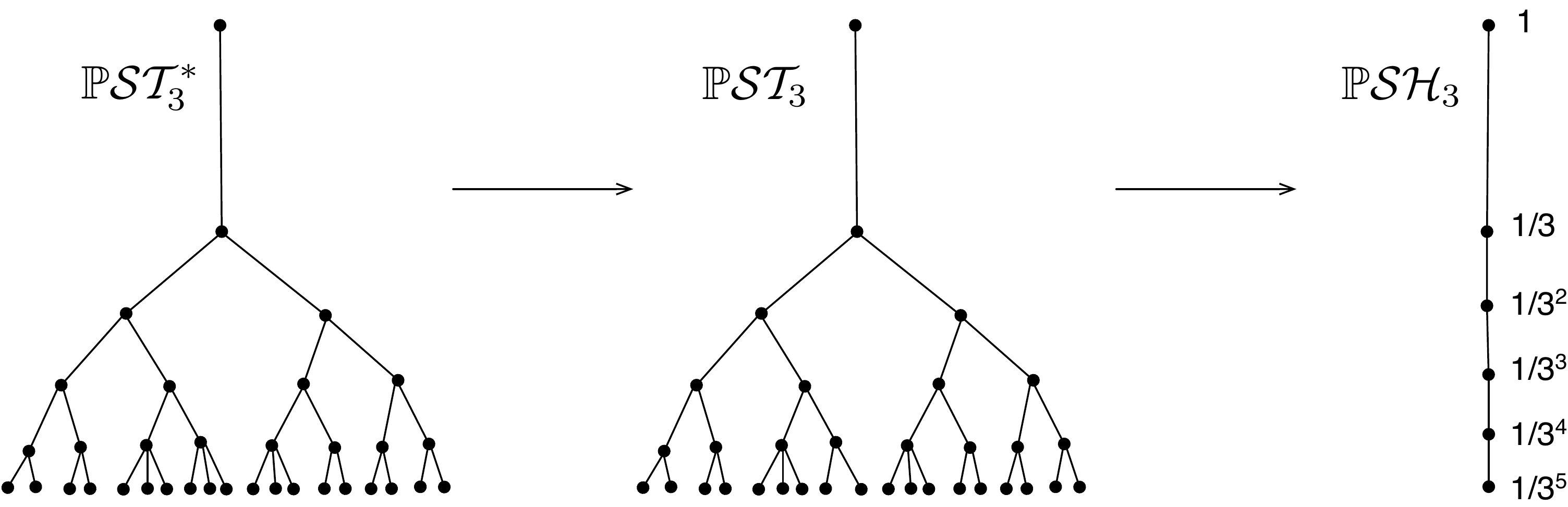}
\caption{The simplicial projections of Theorem \ref{thm:simplicial} in degree $d=3$, drawn for critical heights $1/3^5 \leq h_2/h_1 \leq 1$ (not to scale).  The edges of $\P\cS\cT_3^*$ correspond to structurally stable topological conjugucy classes of polynomials in the shift locus of $\MP_3$.  Note that $\P\cS\cT_3^*$ has one vertex more than $\P\cS\cT_3$ at height $1/3^5$.   }
\label{degree 3 simplicial}
\end{figure}

\begin{theorem}
\label{thm:simplicial}
The spaces $\P\cS\cT_d^*$, $\P\cS\cT_d$, and $\P\cS\cH_d$ carry a  canonical, locally finite simplicial structure, and the projections $\P\cS\cT_d^* \to \P\cS\cT_d \to \P\cS\cH_d$ are simplicial.  
\end{theorem}

By the above-mentioned description of structurally stable maps in the shift locus from \cite{McS:QCIII}, we obtain:

\begin{theorem}
\label{thm:conjugacy_classes}
The set of globally structurally stable topological conjugacy classes of polynomials in the shift locus of $\MP_d$ is in bijective correspondence with the set of 
top-dimensional open simplices in  $\P\cS\cT_d^*$.
 \end{theorem}

It would be useful, therefore, to understand the projectivized shift locus better.  A classification of the stable conjugacy classes in the shift locus is one of the goals of the article \cite{DP:combinatorics}.  In degree $d=2$, the projectivization $\P\cS\cT_2^*$ is a point, while for degree $d=3$ it is an infinite simplicial tree.  See Figure \ref{degree 3 simplicial}.  Using the combinatorics of \cite{DP:combinatorics}, a computational enumeration of the edges of the tree $\P\cS\cT_3^*$ is carried out in \cite{DS:count}.


\bigskip\section{Polynomials}  \label{polynomials}

In this section we define the moduli space $\MP_d$ and present some basic properties of the critical heights map $\cG: \MP_d \to \cH_d$.  We prove Theorem \ref{thm:surjective}, showing that $\cG$ is surjective.

\subsection{Monic and centered polynomials}
Every polynomial 
	$$f(z) = a_0 z^d + a_1 z^{d-1} + \cdots + a_d,$$
with $a_i\in \C$ and $a_0\not=0$, is conjugate by an affine transformation $A(z) = az+b$ to a polynomial which is monic ($a_0=1$) and centered ($a_1 = 0$).  The monic and centered representative is not unique, as the space of such polynomials is invariant under conjugation by $A(z) = \zeta z$ where $\zeta^{d-1} = 1$.  In this way, we obtain a finite branched covering
	$$\cP_d \to \MP_d$$
from the space $\cP_d \iso\C^{d-1}$ of monic and centered polynomials to the moduli space $\MP_d$ of conformal conjugacy classes of polynomials.  Thus, $\MP_d$ has the structure of a complex orbifold of dimension $d-1$.

It is sometimes convenient to work in a space with marked critical points.  Let $H \subset \C^{d-1}$ denote the hyperplane given by $\{c=(c_1, \ldots, c_{d-1}) : c_1 + \ldots + c_{d-1}= 0\}$.  Then the map 
\[ \rho: H \times \C \to \cP_d\]
given by 
\begin{equation}
\label{eqn:param}
 \rho(c_1, \ldots, c_{d-1}; a) = \int_0^z d\cdot \prod_{i=1}^{d-1}(\zeta - c_i)\; d\zeta + a 
 \end{equation}
gives a proper polynomial parameterization of $\cP_d$ by the location of the critical points and the image of the origin.  Setting $\cP^\times_d=H\times \C$, we refer to $\cP^\times_d$ as the space of {\em critically marked} polynomials.  The {\em marked shift locus} is the subset of $\cP_d^\times$ defined by
	$$\cS_d^\times = \{(c\,;a) \in \cP^\times_d: G_f(c_i)>0 \mbox{ for all } i\}$$
where $f = \rho(c_1, \ldots, c_{d-1}; a)$ of equation (\ref{eqn:param}).

\subsection{Critical escape rates} \label{subsec:marked}

The escape rate  $G_f(z) = \lim_{n\to\infty} \frac{1}{d^n} \log^+ |f^n(z)|$ is continuous in $(f,z)\in \cP_d \times \C$, and the {\em maximal escape rate} 
	$$M(f) = \max\{ G_f(c): f'(c)=0\}$$
is a continuous, proper function on $\cP_d$;  see \cite{Branner:Hubbard:1}.  The definition implies that the escape rate satisfies the functional equation $G_f(f^n(z))=d^nG_f(z)$ for all $n\in \N$.  On the open subset of $\cP_d\times\C$ where $G_f(z)>0$, the function $(f,z) \mapsto G_f(z)$ is a locally uniform limit of pluriharmonic functions $d^{-n} \log|f^n(z)|$.  So the map $(f,z)\mapsto G_f(z)$ is pluriharmonic where $G_f(z)>0$.  

Recall that properness of a map is sensitive to changes in the codomain. 

The critical heights map 
	$$\cG: \MP_d \to \cH_d$$
is defined by setting 
	$$\cG(f) = (G_f(c_1), \ldots, G_f(c_{d-1})),$$
where $G_f$ is the escape-rate function of $f$ and $\{c_1, \ldots, c_{d-1}\}$ is the set of critical points of $f$, ordered so that $G_f(c_1) \geq G_f(c_2) \geq \cdots \geq G_f(c_{d-1}) \geq 0$. The critical heights map $\cG$ is proper since it is a factor of the proper map $f \mapsto \cG(f) \mapsto M(f)$.  
We remark that the choice of labeling for the critical points is somewhat artificial; we could just as well have taken $\cG$ to the space of unordered $(d-1)$-tuples.  The decomposition of $\MP_d$ formed by the distinct connected components of fibers of $\cG$ within $\MP_d$ would be the same.  

The critical heights map $\cG$ lifts to a map $\cG^\times: \cP_d^\times \to [0,\infty)^d$ on the space critically marked polynomials by setting
	$$\cG^\times(c\,;a) = (G_f(c_1), \ldots, G_f(c_{d-1})),$$
where $f$ is given by equation (\ref{eqn:param}) and $c=(c_1, \ldots, c_{d-1})$.  The map $\cG^\times$ is continuous since it is a composition of continuous maps.  It is proper since the diagram
$$\xymatrix{ \cP^\times_d \ar[d]_\rho \ar[r]^{\cG^\times} & [0,\infty)^{d-1}\ar[d] \\
		\cP_d   \ar[r]^\cG & \cH_d  }$$
commutes and the left-hand vertical and bottom maps are proper.  We conclude:

\begin{lemma}  \label{G shift}
The restriction of the critical heights map to the marked shift locus
	$$\cG^\times: \cS_d^\times \to (0,\infty)^{d-1}$$
is pluriharmonic and proper.
\end{lemma}


\subsection{Proof of Theorem \ref{thm:surjective}}  It remains only to show surjectivity of $\cG:  \MP_d \to \cH_d$.  Referring to the commutative diagram above, it is enough to show that the lifted map $\cG^\times: \cP_d^\times \to [0,\infty)^{d-1}$ is surjective.  Since it is proper, in order to show that it is surjective, it will suffice to show that the restriction to the marked shift locus $\cG^\times: \cS_d^\times \to (0,\infty)^{d-1}$ is both open and closed.  That it is closed is a consequence of its properness, proved in Lemma \ref{G shift} above.  To prove that it is open, we work with holomorphic maps, and appeal again to Lemma \ref{G shift}.

Given $f \in \cP_d$, the {\em B\"ottcher map} $w=\phi_f(z)$ is the unique analytic isomorphism
\[ \phi_f: \{z : G_f(z) > M(f)\} \to \{ w : |w| > e^{M(f)}\} \]
tangent to the identity near infinity and satisfying $\phi_f\circ f\circ \phi_f^{-1}  (z) = z^d$.  The map $(f, z) \mapsto \phi_f(z)$ is analytic in $f$ and $z$ (see e.g. \cite[Prop 3.7]{Branner:Hubbard:2}).

Fix an integer $n\geq 1$. 
\begin{itemize}
\item Let 
\[ \cS_d^\times(n) = \{ (c;a)\in \cS^\times_d : G_f( f^n(c_i))> M(f), \; 1 \leq i \leq d-1\}\]
where $f=\rho(c;a)$.   That is, $\cS_d^\times(n)$ is the set of critically marked polynomials $f$ for which the $n$th iterates of the critical points of $f$ lie in the domain of the B\"ottcher map $\phi_f$.    By \cite[Lemmas 5.1, 5.2]{DP:basins}, the set $\cS_d^\times(n)$ is connected.  
\item For $w=(w_1, \ldots, w_{d-1}) \in \C^{d-1}$ with each $|w_i| >1$, define $M(w)=\max \log|w_i|$.  For such $w$, set 
\[ \cW^\times_d(n) = \{ w : d^n \log|w_i| > M(w), \; 1 \leq i \leq d-1\}.\]
\item For $h=(h_1, \ldots, h_{d-1}) \in (0,\infty)^{d-1}$, define $M(h) = \max h_i$.  Set 
\[ \cH^\times_d(n) = \{ h : d^n h_i > M(h), \; 1 \leq i \leq d-1\}.\]
Since $\cH^\times_d(n)$ is convex, it is connected.  
\end{itemize}
The functional equation satisfied by the escape rate implies that the function $\Phi_n: \cS_d^{\times}(n) \to \cW_d^\times(n)$ given by 
\[ \Phi_n(c\, ;a) = (\phi_f(f^n(c_1)), \ldots, \phi_f(f^n(c_{d-1})))\]
indeed takes values in $\cW^\times_d(n)$.  Restricted to $\cS_d^\times(n)$, the map $\cG^\times$ is a composition of maps of domains 
\[  \cS_d^\times(n) \stackrel{\Phi_n}{\longrightarrow} \cW_d^\times(n)\stackrel{\Lambda_n}{\longrightarrow} \cH_d^\times(n)\]
where $\Lambda_n(w) =  (\log|w_1|, \ldots, \log|w_{d-1}|)$ is an open map.  
The map $\Phi_n$ is analytic map between domains in $\C^{d-1}$.  It is also proper since it is a factor of the proper map $\cG^\times$; see Lemma \ref{G shift}.   By \cite[Theorems 15.1.15, 15.1.16]{rudin:unitball}, the map $\Phi_n$ is open; we remark that it is also surjective, and that the set of regular values is open.

Since the domains $\cS_d^\times(n)$, $n=1, 2, 3, \ldots$  exhaust $\cS_d^\times$, we conclude that $\cG^\times: \cS_d^\times \to (0,\infty)^{d-1}$ is open.
\qed

\subsection{Remark:  tree-realization implies surjectivity of $\cG$}
One can also prove Theorem \ref{thm:surjective} via the tree realization theorem of \cite{DM:trees}.  Indeed, using the axioms for polynomial trees given there, it is easy to construct an abstract tree $(F, T)$ with any given collection of heights in $\cH_d$.  The simplest construction would place all critical points along a line which joins the Julia set $J(F)$ with $\infty$ and which is fixed by $F$, building the dynamical tree inductively as the height descends.  Because $\cG: \MP_d\to\cH_d$ factors through the space of trees $\cT_d$, we see immediately that $\cG$ is surjective.

\subsection{Critical heights in degree 2}  \label{degree 2 heights}
Finally, we point out that the critical heights map is well-known object in degree $d=2$ where there is a unique critical point.  The moduli space $\MP_2 \iso \C$ is parametrized by the polynomials 
	$$f_c(z) = z^2 + c, \quad c\in\C$$
and the critical heights map $\cG: \C \to \cH_2 = [0,\infty)$ is given by 
	$$\cG(c) = G_c(0).$$
It is known to be equal to $1/2$ times the Green's function for the Mandelbrot set $\cC_2$.  As the Mandelbrot set is connected, all fibers of $\cG$ are connected.  The fibers over positive heights are the equipotential curves, thus each is a simple closed loop around $\cC_2$.  See \cite{Douady:Hubbard}.


\bigskip\section{The quotient space $\cT_d^*$}

Recall that $\cT_d^*$ is the quotient space of $\MP_d$ obtained by collapsing each connected component of a fiber of $\cG: \MP_d \to \cH_d$ to a point.  In this section, we describe some basic topological properties of $\cT_d^*$, and we see that there are only finitely many points in each shift-locus fiber.

\subsection{Degree 2}
As described in \S\ref{degree 2 heights}, the critical heights map in degree 2 is simply (a multiple of) the Green's function for the Mandelbrot set $\cC_2$ in $\MP_2$, so all fibers are connected.  It is immediate to see that the critical heights map induces a homeomorphism
	$$\cT_2^* \longrightarrow \cH_2 =  [0,\infty).$$
Thus in degree 2, the space $\cT_2^*$ coincides with the space of trees $\cT_2$.

\subsection{Monotone-light factor}

\begin{prop} \label{T* is Hausdorff}
The quotient space $\cT_d^*$ is Hausdorff, separable, and metrizable, and the projection $\MP_d \to \cT_d^*$ is proper and closed.
\end{prop}

\proof  Properness follows since the map $\MP_d \to \cT_d^*$ is a factor of the proper map $\MP_d \to \cH_d$.  From \cite[Thm. I.3.5]{daverman:decompositions}, the decomposition of $\MP_d$ into the fibers of $\cG$ is upper semicontinuous, and the quotient space obtained by collapsing fibers of $\cG$ to points is homeomorphic to the image $\cH_d$.  By \cite[Prop. I.4.2]{daverman:decompositions}, the decomposition of $\MP_d$ into the connected components of the fiber of $\cG$ is also upper semicontinuous.   By \cite[Prop. I.1.1,  I.2.1, I.2.2]{daverman:decompositions}, the quotient space $\cT_d^*$ obtained by collapsing these components to points is Hausdorff, separable, and metrizable, and the projection map is closed.
\qed

\subsection{Finiteness of fibers}
By the definition of $\cT_d^*$, the map 
	$\cT_d^*\to \cH_d$ 
has totally disconnected fibers.  We observe here that these fibers are finite when all critical heights are positive.  

\begin{lemma} \label{analytic fibers}  
Fibers of $\cG$ in the shift locus of $\MP_d$ have finitely many connected components.  
\end{lemma}

\proof  
Recall that we can lift $\cG$ to a pluriharmonic map $\cG^\times: \cS_d^\times \to (0,\infty)^{d-1}$ on the critically marked shift locus, as described in \S\ref{subsec:marked}.  The fibers of $\cG^\times$ in $\cS_d^\times$ are locally connected \cite{Bierstone:Milman:semianalytic} and compact (by properness), thus have only finitely many connected components.  
\qed

\medskip
It is well-known that the connectedness locus $\cC_d$ in $\MP_d$ is connected.  See \cite{Douady:Hubbard} for a proof in degree 2, \cite{Branner:Hubbard:1} for degree 3, and \cite{Lavaurs:thesis} for a proof that $\cC_d$ is cell-like in every degree.  Thus, there is a unique point of $\cT_d^*$ with all critical heights equal to 0.  More generally, we have:

\begin{lemma}  \label{uniform heights}
For any height $h\geq 0$, the locus of polynomials in $\MP_d$ with all critical heights equal to $h$ is connected.  Thus, there is a unique point in $\cT_d^*$ with all critical heights equal to $h$.
\end{lemma}

\proof
For $h=0$, the locus is simply the connectedness locus in $\MP_d$, which is well-known to be connected.  For $h>0$, we proved in \cite{DP:basins} that the set $S(f,h)$ is connected for any polynomial $f\in\MP_d$; by definition, this set consists of all maps $g\in\MP_d$ for which the restrictions $f|\{G_f > h\}$ and $g|\{G_g>h\}$ are conformally conjugate and $G_g(c) \geq h$ for all critical points $c$ of $g$.  In particular, for any $f$ with maximal escape rate $M(f)\leq h$, this set $S(f,h)$ is precisely the fiber $\cG^{-1}(h, h, \ldots, h)$ in $\MP_d$.  Proofs for $h>0$ were also given in \cite{Kiwi:combinatorial} and \cite{Dujardin:Favre:critical}.
\qed


\bigskip\section{Trees}  \label{sec:trees}

In this section, we remind the reader about trees and the space of trees. We prove that the critical heights map on $\cT_d$ has totally disconnected fibers, which implies that the projection $\cT_d^*\to \cT_d$ is well-defined.  We also prove Theorem \ref{thm:TstarTfibers} about the projection $\cT_d^*\to \cT_d$.

\subsection{The tree of a polynomial}
Fix a polynomial $f$, and let $G_f: \C\to [0,\infty)$ denote the escape rate, as defined in the Introduction.  The tree associated to $f$ is the monotone-light factor $G_f$:  the quotient $\C\to T(f)$ is obtained by collapsing each connected component of a level set of $G_f$ to a point.  There is a natural simplicial structure on the subset of $T(f)$ which is the quotient of the basin of infinity $X(f)$, realizing this open subset as an infinite locally-finite simplicial tree.  The dynamics of $f$ induces a map $F: T(f) \to T(f)$ which preserves the natural simplicial structure.

\subsection{The space of trees}
We begin with by recalling the topology on the space $\cT_d$ of degree $d$ polynomial trees.  Fix a polynomial $f$ and let $(F, T(f))$ be its tree.  The escape-rate function $G$ of $f$ induces a {\em height function} $h: T(f)\to \R$ and a metric on $T(f)$ so that the distance between adjacent vertices is their difference in height.  Trees $(F_1, T_1)$ and $(F_2, T_2)$ are {\em equivalent} if there is an isometry $i: T_1 \to T_2$ which conjugates $F_1$ to $F_2$.  There is then a unique tree $(F_0, T_0)$ associated to all polynomials with connected Julia set; we call this the {\em trivial tree}.

The topology on $\cT_d$ can be defined as a Gromov-Hausdorff topology, just as on $\mathcal{B}_d$.  The maximal escape rate $M(f)$ of a polynomial is exactly the height of the highest branch point in $T(f)$.  The {\em level} $l\in\Z$ of a vertex $v$ in $T(f)$ is the number of iterates $n$ so that $M(f) \leq h(F^n(v)) < d\, M(f)$.  Roughly speaking, trees are close in the topology on $\cT_d$ if the truncated dynamical systems $(F_N, T_N)$ are close for some large $N$, where $T_N$ is the subtree of vertices of levels $|l| < N$.  

Specifically, a basis of open sets can be defined in terms of two parameters $\eps>0$ and $N\in \N$.  An open neighborhood of a nontrivial tree $(F,T)$ consists of all trees $(F', T')$ for which there is an $\eps$-conjugacy of the restricted trees $(F_N', T_N')$.  Open neighborhoods of the trivial tree $(F_0, T_0)$ are all trees with maximal escape rate $<\eps$.  

With these definitions in place, the proofs of \cite{DM:trees} can be used to show that the projection 
	$ \MP_d \to \cT_d$
is continuous and proper.  See \cite[Theorem 1.4]{DM:trees}, where the projection was only defined on $\MP_d\setminus \cC_d$.

\subsection{Critical heights on trees}

\begin{lemma} \label{THfibers}
The fibers of the critical heights map $\cT_d\to\cH_d$ are totally disconnected; the fibers are finite in the shift locus.
\end{lemma}

\proof

Suppose a fiber $U$ of the critical heights map $\cT_d \to \cH_d$ contains at least two points $(F,T) \neq (F', T')$.  There must be a level $N$ at which the truncated trees $(F_N, T_N), (F'_N, T'_N)$ fail to coincide. There are only finitely many possibilities for these truncated trees.  So $U$ is a disjoint union of relatively open subsets where the truncated tree is constant.  Therefore $U$ is not connected.  

The second statement follows immediately from Lemma \ref{analytic fibers} and the fact that the heights map $\cM_d \to \cH_d$ factors through the space of trees: $\cM_d \to \cT_d \to \cH_d$.  One can also see this directly: by \cite[Theorem 5.7]{DM:trees} a tree in the shift locus is uniquely determined by the subtree down to the level of its lowest critical point.  There are only finitely many combinatorial options once the critical heights are fixed.

\qed

\subsection{Proof of Theorem \ref{thm:TstarTfibers}} 
We first prove that the the map $\cT^*_d \to \cT_d$ is well-defined and that its fibers are totally disconnected. 
Consider the diagram
$$
\xymatrix{
\mathcal{M}_d \ar[rr]^{mon} \ar[dd] \ar[dr]^{mon} &  &\mathcal{T}_d^* \ar[dd]^{light} \ar@{-->}[dl]\\
& \mathcal{M}_d/\sim \ar[dl]_{light} & & \\
\mathcal{T}_d \ar[rr]_{light} & & \mathcal{H}_d 
}
$$
where the middle space is the canonical monotone factor of $\cM_d \to \cT_d$.  Lemma \ref{THfibers} implies that the bottom map is light, so the uniqueness of monotone-light factorizations yields that the spaces $\cM_d/\sim$ and $\cT_d^*$ are the same, i.e. that the map indicated by the dashed arrow is the identity.  


From Lemma \ref{uniform heights}, we know that the fiber of $\cG$ over heights $(h, h, \ldots, h)$ is connected.  Therefore, the fiber of $\cT_d^* \to \cT_d$ over trees with uniform critical heights is a point.  
\qed


\bigskip\section{Quasiconformal deformations}  \label{sec:qc}

In this section, we discuss quasiconformal deformations of polynomials supported on the basin of infinity.  We prove Theorems \ref{generic fibers} and \ref{thm:cone}.  We begin with a continuity statement that we will use several times.





\subsection{Branner-Hubbard motion}  

The upper half-plane $\Hyp=\{\tau = t+is : s>0\}$ may be identified with the subgroup of $\GL_2\R$ consisting of 
matrices
	$${\mtwo{1}{t}{0}{s}}$$ 
with $t\in\R$ and $s>0$, regarded as real linear maps $\tau$ of the complex plane to itself via $\tau \cdot (x+iy) = (x+ty)+i(sy)$.  Note that the parabolic one-parameter subgroup $\{s=1\}$ acts by horizontal shears, while the hyperbolic subgroup $\{t=0\}$ acts by vertical stretches.  The Branner-Hubbard wring motion of \cite{Branner:Hubbard:1} is an action $\Hyp \times \MP_d \to \MP_d$.   

Explicitly, for each polynomial $f$, we consider the holomorphic 1-form $\omega = 2 \, \del G_f$ on the basin $X(f)$.  In the natural Euclidean coordinates of $(X(f), \omega)$, the {\em fundamental annulus}
	$$A(f) = \{M(f) < G_f(z) < d\, M(f)\}$$
may be viewed as a rectangle in the plane, of width $2\pi$ and height $(d-1) M(f)$, with vertical edges identified.  The wringing action is by the linear transformation $\tau$ on this rectangle, transported throughout $X(f)$ by the dynamics of $f$.  For polynomials in the connectedness locus, where $M(f)=0$, the action is trivial.

Put differently, if $\mu$ is the $f$-invariant Beltrami differential $\bar{\omega}/\omega$ on $X(f)$ and 0 elsewhere, we solve the Beltrami equation 
	$$\frac{\delbar \phi_\tau}{\del \phi_\tau} = \frac{-i\tau - 1}{-i\tau + 1} \; \mu$$
for homeomorphism $\phi_\tau: \C\to\C$ and set $f_\tau = \phi_\tau \circ f \circ \phi_\tau^{-1}$.  The map $\tau \mapsto f_\tau$ is analytic in $\tau$, and the escape-rate function of $f_\tau$ satisfies
\begin{equation} \label{stretched G}
	G_{f_\tau}(\phi_\tau(z)) = s \, G_f (z)
\end{equation}
where $s = \Im \tau$.  

The action of the above hyperbolic subgroup is called {\em stretching}, while that of the parabolic subgroup is known as {\em turning}.  Stretching is known to be discontinuous on $\MP_d$ due to the phenomenon of parabolic implosion \cite[Cor. 3.1]{nakane:stretching}; see also \cite{tanlei:discontinuity}, \cite{komori:nakane}.  
Nevertheless we have 

\begin{lemma}
\label{lemma:stretching_continuous}
The Branner-Hubbard wring motion descends to a continuous action on each of the spaces $\mathcal{B}_d$, $\cT_d^*$, $\cT_d$, and $\cH_d$.
In each of these spaces, the point corresponding to the connectedness locus is a global fixed point.  Away from this fixed point the action by stretching is free and proper.   The turning action is nontrivial on $\mathcal{B}_d$ but trivial on $ \cT_d^*$, $\cT_d$, and $\cH_d$.  
\end{lemma}

\proof We begin by verifying that wringing descends continuously to the space $\mathcal{B}_d$.  Suppose $f_n, f $ are polynomials representing elements of $\MP_d$, we fix $\tau_n\to  \tau$ in $\Hyp$, and assume that basins $(X(f_n), f_n) \to (X(f), f)$ in the Gromov-Hausdorff topology on $\mathcal{B}_d$.  Since the forgetful map $\MP_d \to \mathcal{B}_d$ is proper, we may assume $f_n \to f$ for concreteness.  Let $\mu_n$  be the $f_n$-invariant Beltrami differentials on $\C$ corresponding to the wring of $X(f_n)$ by $\tau_n$; they are supported on $X(f_n)$.  Similarly, let $\mu$ correspond to $\tau$; it is $f$-invariant and supported on $X(f)$.   Suitably normalized, there are  unique quasiconformal maps $h_n: \C \to \C$ and $h: \C \to \C$ whose complex dilations are almost everywhere $\mu_n$ and $\mu$, respectively \cite[Prop. 6.1]{Branner:Hubbard:1}.  Let $g_n = h_n\circ f_n \circ h_n^{-1}$ and $g = h\circ f \circ h^{-1}$.  By construction, the maps $h_n$ are uniformly quasiconformal, so we may pass to a subsequence so that $h_n \to \widetilde{h}$ and $g_n \to \widetilde{g}$.   It follows that $h$ conjugates $f$ to $g$, while $\widetilde{h}$ conjugates $f$ to $\widetilde{g}$.   Let $\widetilde{\mu}$ be the complex dilatation of $\widetilde{h}$.  

By \cite[Lemma 10.1]{Branner:Hubbard:1}, $\mu=\widetilde{\mu}$ on $X(f)$.  It follows $g$ and $\widetilde{g}$ are holomorphically conjugate on their basins of infinity, thus $g = \widetilde{g}$ in $\mathcal{B}_d$.  So wringing is continuous on $\mathcal{B}_d$.

By equation (\ref{stretched G}), wringing scales the critical height vector associated to each polynomial.  It follows that on $\mathcal{B}_d$, wringing preserves the fibers of $\mathcal{B}_d \to \cT_d^*$ and so descends continuously to an action on $\cT_d^*$.  The remaining assertions follow easily.  
\qed

\subsection{Multitwisting and the Teichm\"uller space}  \label{twist}
Given $f \in \MP_d$, there is a canonical space of marked quasiconformal deformations of $f$ supported on the basin of infinity.  The general theory, developed in \cite{McS:QCIII}, shows that this space admits the following description.  

Fix a polynomial representative $f: \C\to\C$ of its conjugacy class. The {\em foliated equivalence class} of a point $z$ in the basin $X(f)$ is the closure of its grand orbit $\{w \in X(f): \exists \; n, m\in \Z, f^n(w) = f^m(z)\}$ in $X(f)$.  Let $N$ be the number of distinct foliated equivalence classes containing critical points of $f$.   These critical foliated equivalence classes subdivide the fundamental annulus $A(f)$ into $N$ {\em fundamental subannuli} $A_1, \ldots, A_N$ linearly ordered by increasing height.   It turns out one can define wring motions via affine maps on each of the subannuli $A_j$ independently so that the resulting deformation of the basin $X(f)$ is continuous and well-defined.   The deformations of each subannulus are parameterized by $\Hyp$, so we obtain an analytic map 
\[ \Hyp^N \to \MP_d. \]
By varying the map $f$ as well, we get an action 
\[ \Hyp^N \times \MP_d^N \to \MP_d^N\]
where now $\MP_d^N$ is the locus of maps with exactly $N$ critical foliated equivalence classes.    The Branner-Hubbard wring by $\tau = t + i s \in\Hyp$ applied to $f\in \MP_d^N$ is the action of vector 
$$\left( \frac{2\pi \, m_1\, t}{(d-1) M(f)} + is \, , \, \ldots \, , \, \frac{2\pi \, m_N\, t}{(d-1) M(f)} + is \right) \in\Hyp^N$$ 
where $m_j$ is the modulus of the $j$-th subannulus of $A(f)$, so that 
	$$\sum_j m_j = (d-1) M(f) / 2\pi.$$
The action of $\R^N$ by the parabolic subgroup in each factor is called {\em multitwisting}.   By construction, the multitwisting deformations preserve critical heights.  

Let $\cS_d^N$ be the locus consisting of polynomials in the shift locus with exactly $N$ distinct critical foliated equivalences classes.  

\begin{lemma}
\label{lemma:mult_stretching_continuous}
The wring motion on each fundamental subannulus defines a continuous action 
\[ \Hyp^N \times  \cS_d^N \to \cS_d^N\]
on the subset of the shift locus consisting of polynomials with exactly $N$ distinct critical foliated equivalences classes.  For each $f$, the orbit map $\Hyp^N\times\{f\} \to \cS_d^N \subset \MP_d$ is analytic, and 
the stabilizer of $f$ contains a lattice of translations in $\R^N$ in general and equals a lattice of translations when $N = d-1$.
\end{lemma}

\proof Continuity follows from the same arguments given in the proof of Lemma \ref{lemma:stretching_continuous}; analyticity follows from the analytic dependence of the solution of the Beltrami equation.   To prove the second assertion, suppose $f \in \cS_d^N$.  Let $A_j$ be the $j$th fundamental subannulus of $f$.  Since $f$ belongs to the shift locus, 
\[ d_j = \mbox{\rm lcm}\{ \mbox{\rm deg}(f^n: \widetilde{A} \to A_j) \} < \infty \]
where the least common multiple is taken over all connected components of all iterated inverse images $\widetilde{A}_j$ of $A_j$.  Let $h_j: A_j \to A_j$ be the $d_j$-th power of a right Dehn twist which is affine in the natural Euclidean coordinates on $A_j$.   From the definition of $d_j$ it follows that $h_j$ extends uniquely to a quasiconformal deformation $h_j: X(f) \to X(f)$ which commutes with $f$.  It follows that  the stabilizer of $f$ contains the lattice in $\R^N$ generated by 
	$$(0, \ldots, 0, d_j/m_j, 0, \ldots, 0),$$
where $m_j$ is the modulus of $A_j$, $j=1, \ldots, N$.  

Conversely, any power of a Dehn twist in $A_j$ that stabilizes $f$ arises in this way.  If a multitwist stabilizes $f$ and in some $A_j$ is not a power of a Dehn twist, then it must permute nontrivially a set of marked points on $\partial A_j$ corresponding to the orbits of critical points; this is impossible if the heights of $f$ are independent.
\qed

\subsection{Proof of Theorem \ref{generic fibers}}  \label{generic proof}
We now prove that for generic critical heights in the shift locus, the fibers of $\cG$ are tori coinciding with twist orbits.  In other words, the quotient space $\cT_d^*$ is generically the orbit space for the twist deformation.  We begin by characterizing these generic critical heights.

A critical height value $(h_1, \dots, h_{d-1})$ is {\em generic} for $\cG$ if
\begin{enumerate}
\item	$h_i >0$ for all $i$, and
\item	$h_i \not= d^n h_j$ for each $i\not= j$ and all $n\in \Z$.
\end{enumerate}
It is clear that condition (1) corresponds to the shift locus, while condition (2) corresponds to the stratum $\cS_d^{d-1}$ where all critical points lie in distinct foliated equivalence classes.  In this stratum, the twisting deformation space has maximal dimension $d-1$.  Lemma \ref{lemma:mult_stretching_continuous} implies that $f$ lies in a $(d-1)$-dimensional analytic torus of maps obtained from multitwisting.

On the other hand, the existence of $d-1$ independent stretches in the fundamental subannuli implies that $\cG$ has maximal rank $d-1$ along a generic fiber:  

\begin{lemma} \label{Grank}
For each $f\in \cS_d^N$, we have $\rank_f \mathcal{G} \geq N$.
\end{lemma}

\proof The heights map $\cG: \cS_d \to \cS\cH_d$ is smooth.  The locus $\cS\cH_d^N \subset \cH_d$ is a smooth submanifold of dimension $N$.  By Lemma \ref{lemma:mult_stretching_continuous}, for each $f \in \cS_d^N$, independent stretching of the fundamental subannuli of $f$ provides an analytic section of $\cG$ over a neighborhood of $\cG(f)$ in $\cS\cH_d^N$ and so $\rank_f\cG \geq N$.  

In detail, let $G$ be the escape rate of $f$.  For each escaping critical point $c$, let $l(c)$ be its level:  the least integer such that $G(f^{l(c)}(c))\geq G(c_1)$.  Choose a critical point $c_j$ representing each foliated equivalence class, and relabel them as $c_1, \ldots, c_N$ so that 
	$$G(f^{l(c_j)}(c_j)) < G(f^{l(c_{j+1})}(c_{j+1}))$$
for all $j< N$.  The moduli $m_j$ of the fundamental subannuli $A_j$ satisfy
	$$m_1 = \frac{d-1}{2\pi} \;G(c_1)$$
when $N=1$; otherwise,
\begin{equation} \label{moduli matrix}
	\left( \begin{array}{c} m_1 \\ m_2 \\ m_3 \\ \vdots \\ m_N \end{array} \right)
		= \frac{1}{2\pi} 
		\left(  \begin{array}{ccccc} 
		-1 & 1 & 0 &  & 0 \\
		0 & -1 & 1 & \cdots & 0 \\
		0 & 0 & -1 &  &  0 \\
		 & \vdots & & \ddots & \\
		d & 0 & 0 &  & -1 \end{array} \right) 
		\left( \begin{array}{c} 
		G(c_1) \\ G(f^{l(c_2)}(c_2)) \\ G(f^{l(c_3)}(c_3)) \\ \vdots \\ G(f^{l(c_N)}(c_N))
		\end{array} \right)
\end{equation}
The matrix has determinant $(-1)^{N-1}(d-1)$, so it is invertible.  The moduli $m_j$ can be freely adjusted under analytic stretching deformations.  Recalling that $G(f^k(z)) = d^k\, G(z)$, we see that the rank of $\cG$ is at least $N$ at $f$.
\qed

\medskip
As a consequence of Lemma \ref{Grank}, the generic fiber of $\cG$ is a smooth compact submanifold of real dimension $d-1$.    Because the twist orbit is both open and closed in the fiber, it must coincide with a connected component.

We complete the proof of Theorem \ref{generic fibers} by observing that the conclusion only holds for the generic heights defined above.  Indeed, for non-generic critical heights, the dimension of the twisting deformation is $< d-1$.  
\qed

\subsection{Proof of Theorem \ref{thm:cone}}
The existence of the cone structures on spaces $\mathcal{B}_d$, $\cT_d^*$, $\cT_d$, and $\cH_d$, and their basic properties follow immediately from Lemma \ref{lemma:stretching_continuous} and equation (\ref{stretched G}).  
\qed


\bigskip\section{The projectivization $\P\cT_d^*$}

Recall from Theorem \ref{thm:cone} that there is a well-defined projectivization $\P\cT_d^*$ of the space $\cT_d^*$ via stretching.  
In this section we prove Theorem \ref{thm:contractible} which states that $\P\cT_d^*$ is compact and contractible.  The proof is very similar to the proof in \cite{DM:trees} that the projectivized space of trees $\P\cT_d$ is compact and contractible.  Compactness will follow easily from the properness of $\cG$.  To prove contractibility, the idea is to construct a deformation retract of $\P\cT_d^*$ which ``lifts" the corresponding retract of $\P\cT_d$.  

\subsection{The root point of $\P\cT_d^*$}
Recall that the projectivization $\P\cT_d^*$ is identified with the slice $\cT_d^*(1)$ of $\cT_d$.  From Lemma \ref{uniform heights}, we know that there is a unique point in $X(h) \in \cT_d^*$ corresonding to maps with uniform critical heights $=h$, for any $h\geq 0$.  Let $X(1)\in\cT_d^*(1)$ be this point for height $h=1$.  We will call this the {\em root point} of $\P\cT_d^*$.

\subsection{Paths through the shift locus} \label{subsection:paths}
Assume $f\in\MP_d$ lies in the shift locus and has maximal escape-rate $M(f)=1$.  We begin by recalling the construction of \cite[Lemma 5.2]{DP:basins} which defines a path $f_t$ from $f$ to a polynomial with {\em all} critical heights equal to 1, in such a way that  for every time $0 \leq t \leq 1$, the restriction $f|\{G_f > t\}$ is conformally conjugate to $f_t|\{G_{f_t} > t\}$ and $f_t$ satisfies 
	$$G_{f_t}(c) \geq t$$
for all critical points $c$ of $f_t$.  The idea is to push the lowest critical values of $f$ up along their external rays until they reach height $d$.  This pushing deformation is uniquely defined as long as the critical values do not encounter the critical points of $G_f$.  When this occurs, multiple external rays land at a critical value, so there are finitely many choices for continuing the path.  Nevertheless, for any choice and every $t > 0$, the polynomial $f_t$ lies in a {\em connected} set we called $S(f,t)$; the set $S(f,t)$ consists of {\em all} polynomials $g\in \MP_d$ which are conformally conjugate to $f$ on the restriction $\{G_f > t\}$ and have all critical heights $\geq t$.  The set $S(f,t)$ is thus contained in a single connected component of a fiber of the critical heights map $\cG$.  

We will use this pushing-up deformation to define {\em canonical} paths in $\P\cT_d^*$.  

\begin{lemma} \label{converging to S(f,t)}
Suppose $f_n$ is a sequence of polynomials converging to $f$ as basins in $\mathcal{B}_d$, and suppose $t_n$ is a sequence of real numbers converging to $t>0$.  Then every subsequential limit in $\MP_d$ of polynomials $g_n \in S(f_n,t_n)$  must lie in $S(f,t)$.  
\end{lemma}

\proof
By hypothesis, the polynomial $g_n$ is conformally conjugate to $f_n|\{G_{f_n}> t_n\}$ on its restriction to $\{G_{g_n} > t_n\}$.  It also satisfies $G_{g_n}(c) \geq t_n$ at all critical points $c$.  The convergence of $f_n \to f$ in $\mathcal{B}_d$ means that the restrictions $f_n|\{G_{f_n} > 0\}$ converge to $f|\{G_f > 0\}$ in the Gromov-Hausdorff sense.  See \cite{DP:basins} for details.  In particular, the restrictions $f_n|\{G_n > t_n\}$ must converge to $f|\{G_f > t\}$ (geometrically).  It follows that any uniform limit $g$ of the polynomials $g_n$ (as elements of $\MP_d$) has restriction $g|\{G_g > t\}$ which is conformally conjugate to $f|\{G_f > t\}$.   Furthermore, the critical heights of $g_n$ converge to those of $g$, so we can conclude that $g \in S(f,t)$.   
\qed

\begin{lemma} \label{Q_t}
Fix an element $q \in \cT_d^*$, and let $Q$ be the fiber over $q$ in $\MP_d$.  Then for all $t>0$, the union 
	$$Q_t = \bigcup_{f\in Q} S(f,t)$$
is connected and lies in a unique connected component of a fiber of the critical heights map $\cG$.  
\end{lemma}

\proof
By construction, the fiber $Q$ is connected.  Also from the definitions, the critical heights map $\cG$ is constant on the set $Q_t$.  It suffices to show that $Q_t$ is connected.  Recall that for each fixed $f$, the set $S(f,t)$ is connected.  Let $A$ and $B$ be disjoint open sets such that $Q_t \subset A\cup B$.  If $Q_t$ has nonempty intersection with both $A$ and $B$,  then at least one of the sets 
	$$Q_A = \{f \in Q: S(f,t) \subset A\}$$
or
	$$Q_B = \{f \in Q: S(f,t) \subset B\}$$
is not closed.  Suppose it is $Q_A$.  Since $Q = Q_A \cup Q_B$ is compact and connected, there is a sequence $f_n\in Q_A$ which converges to $f \in Q_B$.   By Lemma \ref{converging to S(f,t)}, elements of $S(f_n, t) \subset A$ must converge to $S(f,t) \subset B$.  This is a contradiction, so $Q_t$ must be connected.
\qed

\subsection{Proof of Theorem \ref{thm:contractible}}
Compactness follows immediately from the properness of $\cG$.
Indeed, $\P\cT_d^*$ is homeomorphic to the slice $\cT_d^*(1)$ of $\cT_d^*$ consisting of points with maximal critical height equal to 1.  The slice $\cT_d^*(1)$ is the quotient of the preimage $$\cG^{-1}(1, h_2, \ldots, h_{d-1})$$ in $\MP_d$ where $1 \geq h_2 \geq \cdots \geq h_{d-1} \geq 0$.  By properness of $\cG$, this locus is compact, so its quotient in $\cT_d^*$ is also compact.

For contractibility, we define a map 
	$$R^*: [0,1] \times \P\cT_d^* \to \P\cT_d^*$$
by setting $R^*(0,q) = q$ for all $q\in \cT_d^*(1)$ and $R^*(t, q) = q_t$ for $t>0$, where $q_t$ is the unique element in $\cT_d^*(1)$ associated to the set $Q_t$ of Lemma \ref{Q_t}.  We observed in the proof of \cite[Proposition 4.6]{DP:basins} that any family of polynomials $f_t \in S(f,t)$ converges to $f$ as $t\to 0$ in the Gromov-Hausdorff topology on $\mathcal{B}_d$.  Consequently, the elements $q_t$ converge to $q$ in the quotient space $\P\cT_d^*$.  Further, the convergence is easily seen to be uniform on a compact neighborhood of $q$, by the structure of the open sets which define the topology of $\cB_d$.  It follows that $R^*$ is continuous at $t=0$.  

It remains to show continuity of $R^*$ for all $t>0$.  Suppose $q_n$ is a sequence in $\cT_d^*(1)$ converging to $q\in \cT_d^*(1)$, and choose any sequence $t_n \to t$ in $(0,1]$.   For any elements $f_n\in \MP_d$ in the fiber $Q_n$ over $q_n$, pass to a subsequence so that $f_n \to f$ in $\mathcal{B}_d$.  By continuity, $f$ must project to $q$ in $\cT_d^*(1)$.  By Lemma \ref{converging to S(f,t)}, any sequence $g_n \in S(f_n, t_n)$ has subsequential limits inside $S(f,t)$.  Note that $g_n$ projects to $(q_n)_{t_n}$ and all of $S(f,t)$ projects to $q_t$.  Therefore, $R^*(t_n,q_n) \to R^*(t, q)$, so $R^*$ is continuous.
\qed



\bigskip\section{Simplicial structures and stable conjugacy classes} 

In this section, we prove Theorems \ref{thm:simplicial} and \ref{thm:conjugacy_classes}.  In particular, we define
a locally finite simplicial complex structure on $\P\cS\cT_d^*$ and show that the open simplices of top dimension correspond to the topological conjugacy classes of structurally stable polynomials in the shift locus.  Our description of the simplicial structure differs from the proof in \cite{DM:trees} of the analogous statement for $\P\cS\cT_d$, since we do not have an intrinsic description of the points of $\P\cS\cT_d^*$.

\subsection{The open simplices}  \label{open simplices}
For each $N \in\{1, \ldots, d-1\}$, let 
	$$\sigma_N = \{(x_1, \ldots, x_N)\in\R^N:  x_i > 0\ \forall i, \  x_1 + \cdots + x_N = 1 \}$$
be an {\em open} simplex of dimension $N-1$.  Fix $f$ in the shift locus $\cS_d$ with exactly $N$ critical foliated equivalence classes and maximal critical height $M(f) = 1$.   As described in \S\ref{twist}, the fundamental subannuli $A_1, \ldots, A_N$ of $f$ can be independently stretched, to freely adjust the moduli of the $A_j$.  We define a continuous, injective map 
	$$\sigma_N \to \cS_d$$
sending $(x_1, \ldots, x_N)$ to the unique point in the stretch-orbit of $f$ with 
	$$\mod A_j = \frac{(d-1)\ x_j}{2\pi}.$$  
The image of $\sigma_N$ contains $f$ and lies in the stratum $\cS_d^N(1)$ consisting of polynomials in the shift locus with exactly $N$ fundamental subannuli and maximal critical height = 1.  

The critical heights of all maps in the image of $\sigma_N$ can be computed directly from $\cG(f)$ and $(x_1, \ldots, x_N)$, using equation (\ref{moduli matrix}).  It follows that the composition 
	$$\sigma_N \to \cS_d(1) \to \cT_d^*(1) \to \cT_d(1) \to \cH_d(1)$$
is injective.  

In fact, the image of $\sigma_N$ is an entire {\em connected component} of the stratum $\P\cS\cH_d^N$.  We shall see that these $\sigma_N$ fit together to define the simplicial structure on each of the projectivizations $\P\cS\cT_d^* \to \P\cS\cT_d \to \P\cS\cH_d$, and the projection maps are simplicial.  We begin by describing the simplicial structure on $\P\cS\cH_d$ in detail.

\subsection{Simplicial structure on the space of critical heights} \label{heights structure}
Recall that the projectivization $\P\cS\cH_d$ is naturally identified with the set of vectors
	$$\{h = (1, h_2, \ldots, h_{d-1})\in \R^{d-1}:  0 < h_{d-1} \leq h_{d-2} \leq \cdots \leq h_2 \leq 1\}.$$
In \S\ref{strat} of the Introduction, we introduced a partition 
	$$\P\cS\cH_d = \P\cS\cH_d^1 \ \sqcup \cdots \sqcup \ \P\cS\cH_d^{d-1}$$ 
based on the maximal number of {\em independent} heights in degree $d$.  Positive numbers $x$ and $y$ are independent if $x \not= d^n y$ for any integer $n$.

\begin{figure} 
\includegraphics[width=3.0in]{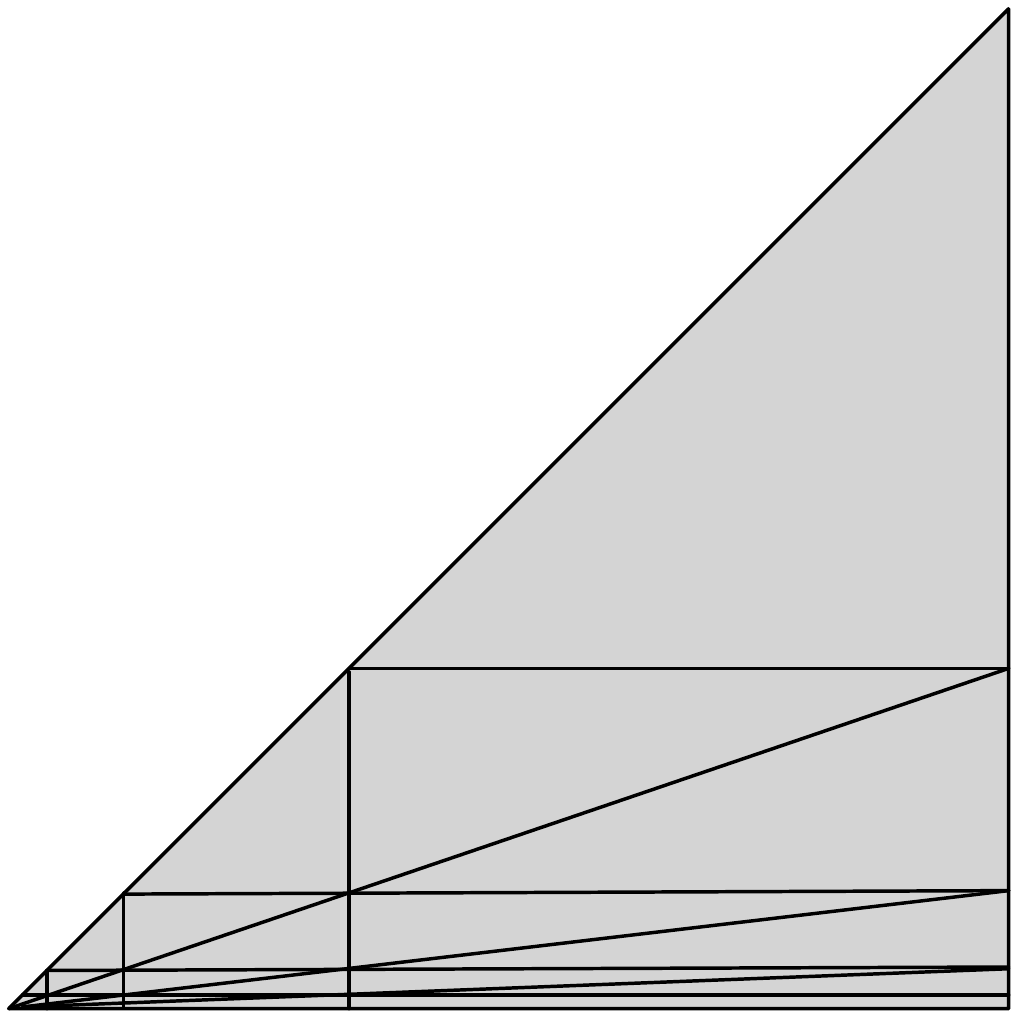}
\caption{degree 4 simplicial structure on $\P\cS\cH_4$, depicted as the triangle $\{0 < y \leq x \leq 1\}$ in $\R^2$.  Though not drawn to scale, the lines represent height relations $\{x = 1/4^n\}$, $\{y = 1/4^n\}$, and $\{y = x/4^n\}$ for all positive integers $n$. }
\label{degree 4 simplicial}
\end{figure}

The height vectors in $\P\cS\cH_d^1$ form a discrete set and will be the 0-skeleton of $\P\cS\cH_d$.  Now let $h = (1, h_2, \ldots, h_{d-1})$ be a vector with exactly two independent heights.  Choose $h_j$ independent from $1$, and let $n_j$ be the unique integer such that 
	$$1 < d^{n_j}h_j < d.$$
There is a unique continuous map
	$$(0,1) \to \P\cS\cH_d$$
with image containing $h$, sending $x\in (0,1)$ to a height vector with exactly two independent heights and $j$-th coordinate equal to
	$$h_j = \frac{1 + (d -1)\ x}{d^{n_j}} \, .$$
Thus, the map is a bijection from the open 1-simplex $(0,1)$ to the component of $\P\cS\cH_d^2$ containing $h$; it is easy to see that it extends continuously to a simplicial map
	$$[0,1] \to \P\cS\cH_d^1 \sqcup \P\cS\cH_d^2.$$
In this way, we define the locally-finite simplicial structure inductively on strata, so that the union $\P\cS\cH_d^1 \sqcup \cdots \sqcup \P\cS\cH_d^N$ forms the $(N-1)$-skeleton of $\P\cS\cH_d$.  The simplices are shown for degree 3 in Figure \ref{degree 3 simplicial} and for degree 4 in Figure \ref{degree 4 simplicial}.

\subsection{Proof of Theorem \ref{thm:simplicial}}
Fix a polynomial $f$ in the shift locus $\cS_d$, with maximal escape rate $M(f) = 1$ and $N$ critical foliated equivalence classes.  Let 
	$$\overline{\sigma_N} = \{(x_1, \ldots, x_N)\in\R^N:  x_i \geq 0\ \forall i, \  x_1 + \cdots + x_N = 1 \}$$
be a closed simplex of dimension $N-1$.  We map the interior $\sigma_N$ to the stretching orbit of $f$ within $\cS_d(1)$, as described above in \S\ref{open simplices}.  

We have already seen in \S\ref{heights structure} that the maps $\sigma_N\to \P\cS\cH_d^N$ extend to continuous, injective maps  
	$$\overline{\sigma_N} \to\P\cS\cH_d$$ 
to define a simplicial structure on $\P\cS\cH_d$.  Because the fibers of $\cT_d^* \to \cT_d\to\cH_d$ are finite in the shift locus, by Lemma \ref{analytic fibers}, the maps from $\sigma_N$  extend continuously and injectively through the sequence of spaces
	$$\overline{\sigma_N} \to \P\cS\cT_d^* \to \P\cS\cT_d \to \P\cS\cH_d.$$

It follows from \cite[\S2]{DM:trees}, and specifically the Proposition 2.17 there, that the $(N-1)$-dimensional open simplices $\sigma_N \to \P\cT_d$ either coincide or are disjoint, because the height metrics on a given combinatorial tree are parameterized by $\sigma_N$.  Using the fact that the extended maps $\overline{\sigma_N}\to\P\cS\cH_d$ define the simplicial structure on the space of heights, we conclude that the closed simplices $\overline{\sigma_N} \to \P\cT_d$ provide the maps for a well-defined simplicial structure on $\P\cS\cT_d$, compatible with the projection $\P\cS\cT_d \to \P\cS\cH_d$.  See \cite[Lemma I.2.1]{munkres:algebraictopology}.  See also \cite[\S11]{DM:trees} where the simplicial structure on $\P\cS\cT_3$ is described in detail.

It remains to show that the images of two open simplices in $\P\cS\cT_d^*$ either coincide or are disjoint.  Fix a critical height vector $h$ in $\P\cS\cH_d$, and let $Q \subset \cS_d(1)$ be a connected component of a fiber of $\cG$ over $h$.  The continuity of stretching from Lemma \ref{lemma:stretching_continuous} implies that the image of $Q$ under any stretch must remain in a connected component of a fiber of $\cG$.  This shows immediately that the image of simplices either coincide or are disjoint in $\P\cS\cT_d^*$.  
\qed

\subsection{Remark: closed simplices in $\cS_d$}
While not needed for the proof of Theorem \ref{thm:simplicial}, it is useful to observe that the open simplices $\sigma_N \to \cS_d$
defined in \S\ref{open simplices} extend continuously and injectively to the closure 
	$$\overline{\sigma_N} \to \cS_d$$
in the space of polynomials, not only in the quotient spaces.  One way to see this is to use the Gromov-Hausdorff topology on the restrictions $(f, X(f))$ to the basins of infinity; this topology is equivalent to the topology of uniform convergence on $\MP_d$ in the shift locus \cite{DP:basins}.  The stretching deformation adjusts the relative heights of annuli in the metric on $X(f)$, leaving other invariants of the dynamics of $f$ unchanged.  In particular, the external angles of rays landing at critical points are fixed, and the collection of polynomials with given critical heights and external angles form a discrete set in the shift locus.  It follows immediately that the extension to $\overline{\sigma_N} \to \cS_d$ is well-defined and continuous.

\subsection{Proof of Theorem \ref{thm:conjugacy_classes} }
We conclude this article with the proof that the set of globally structurally stable topological  conjugacy classes of maps in the shift locus is in bijective correspondence with the set of  top-dimensional open simplices of $\P\cS\cT_d^*$.

Recall that for maps in the shift locus, quasiconformal and topological conjugacy classes coincide and are connected.  Further, the structurally stable polynomials coincide with the top-dimensional stratum $\cS_d^{d-1}$ \cite{McS:QCIII}.  Let $f$ be a map with $d-1$ independent heights in the shift locus and let $\sigma_f$ be the open simplex containing the image of $f$ in $\P\cS\cT_d^*$; thus $\sigma_f$ is an open simplex of maximal dimension.  The connectedness of topological conjugacy classes and the definition of the simplicial structure on $\P\cS\cT_d^*$ shows that $f \mapsto \sigma_f$  is well-defined on topological conjugacy classes $[f]$.  It is clearly surjective.    It remains to show injectivity.  Suppose $\sigma_{f_1}=\sigma_{f_2}$.  Then there are elements in the topological conjugacy classes of $f_1$ and $f_2$ which lie in the same connected component of a fiber of $\cG$.  For generic heights, these connected components coincide with twist-deformation orbits by Theorem \ref{generic fibers}, and therefore $[f_1] = [f_2]$.  
\qed


\bigskip\bigskip
\def\cprime{$'$}


\begin{thebibliography}{BDK}

\bibitem[BM]{Bierstone:Milman:semianalytic}
E.~Bierstone and P.~D. Milman.
\newblock {Semianalytic and subanalytic sets}.
\newblock {\em Inst. Hautes \'Etudes Sci. Publ. Math.} {\bf 67}(1988), 5--42.

\bibitem[BDK]{Blanchard:Devaney:Keen}
P.~Blanchard, R.~L. Devaney, and L.~Keen.
\newblock {The dynamics of complex polynomials and automorphisms of the shift}.
\newblock {\em Invent. Math.} {\bf 104}(1991), 545--580.

\bibitem[BH1]{Branner:Hubbard:1}
B.~Branner and J.~H. Hubbard.
\newblock {The iteration of cubic polynomials. {I}. {T}he global topology of
  parameter space}.
\newblock {\em Acta Math.} {\bf 160}(1988), 143--206.

\bibitem[BH2]{Branner:Hubbard:2}
B.~Branner and J.~H. Hubbard.
\newblock {The iteration of cubic polynomials. {II}. {P}atterns and
  parapatterns}.
\newblock {\em Acta Math.} {\bf 169}(1992), 229--325.

\bibitem[Da]{daverman:decompositions}
R.~J. Daverman.
\newblock {\em Decompositions of manifolds}.
\newblock AMS Chelsea Publishing, Providence, RI, 2007.
\newblock Reprint of the 1986 original.

\bibitem[DM]{DM:trees}
L.~DeMarco and C.~McMullen.
\newblock {Trees and the dynamics of polynomials}.
\newblock {\em Ann. Sci. {\'E}cole Norm. Sup.} {\bf 41}(2008), 337--383.

\bibitem[DP1]{DP:basins}
L.~DeMarco and K.~Pilgrim.
\newblock {Polynomial basins of infinity}.
\newblock {\em {\em Preprint, 2009}}.

\bibitem[DP2]{DP:hausdorff}
L.~DeMarco and K.~Pilgrim.
\newblock {Hausdorffization and polynomial twists}.
\newblock {Submitted for publication, 2009}.

\bibitem[DP3]{DP:combinatorics}
L.~DeMarco and K.~Pilgrim.
\newblock {Escape combinatorics for polynomial dynamics}.
\newblock {\em {\em Preprint, 2009}}.

\bibitem[DS]{DS:count}
L.~DeMarco and A.~Schiff.
\newblock {Enumerating the basins of infinity for cubic polynomials}.
\newblock {\em {\em Submitted for publication, 2008}}.

\bibitem[DH]{Douady:Hubbard}
A.~Douady and J.~H. Hubbard.
\newblock {It\'eration des polyn\^omes quadratiques complexes}.
\newblock {\em C. R. Acad. Sci. Paris S\'er. I Math.} {\bf 294}(1982),
  123--126.

\bibitem[DF]{Dujardin:Favre:critical}
R.~Dujardin and C.~Favre.
\newblock {Distribution of rational maps with a preperiodic critical point}.
\newblock {\em Amer. J. Math.} {\bf 130}(2008), 979--1032.

\bibitem[Ki]{Kiwi:combinatorial}
J.~Kiwi.
\newblock {Combinatorial continuity in complex polynomial dynamics}.
\newblock {\em Proc. London Math. Soc. (3)} {\bf 91}(2005), 215--248.

\bibitem[KN]{komori:nakane}
Y.~Komori and S.~Nakane.
\newblock {Landing property of stretching rays for real cubic polynomials}.
\newblock {\em Conform. Geom. Dyn.} {\bf 8}(2004), 87--114 (electronic).

\bibitem[La]{Lavaurs:thesis}
P.~Lavaurs.
\newblock {\em Systemes dynamiques holomorphes: explosion de points periodiques
  paraboliques}.
\newblock Thesis, Orsay, 1989.

\bibitem[McS]{McS:QCIII}
C.~T.~McMullen and D.~P.~Sullivan.
\newblock {Quasiconformal homeomorphisms and dynamics. {I}{I}{I}. {T}he
  {T}eichm\"uller space of a holomorphic dynamical system}.
\newblock {\em Adv. Math.} {\bf 135}(1998), 351--395.

\bibitem[Mu]{munkres:algebraictopology}
J.~R.~Munkres.
\newblock {\em Elements of algebraic topology}.
\newblock {Addison-Wesley Publishing Company, Menlo Park, CA, 1984}.

\bibitem[Na]{nakane:stretching}
S.~Nakane.
\newblock {Branner-{H}ubbard-{L}avaurs deformations for real cubic polynomials
  with a parabolic fixed point}.
\newblock {\em Conform. Geom. Dyn.} {\bf 13}(2009), 110--123.

\bibitem[R]{rudin:unitball}
W.~Rudin.
\newblock{Function theory in the unit ball of {$\mathbb C\sp n$}}.
\newblock{Springer-Verlag, Berlin, 2008}.
\newblock{Reprint of the 1980 edition}.

\bibitem[T]{tanlei:discontinuity}
Tan~Lei.
\newblock{Stretching rays and their accumulations, following Pia Willumsen.}
\newblock{ In {\em Dynamics on the Riemann sphere}, 183-208}.
\newblock{ Eur. Math. Soc., Z\"urich } 2008.

\end{thebibliography}

 \end{document}